%% file: remarks-vanishing.tex
\documentclass[a4paper,12pt]{amsart}
\usepackage{amssymb}
\usepackage{amscd}
\usepackage[abbrev,alphabetic]{amsrefs}
\usepackage[hyperfootnotes=false]{hyperref}
	\hypersetup{
	hidelinks,
	pdfauthor={Fabio Bernasconi},
	pdftitle={vanishing-char-p},
}
\usepackage{xypic}
\usepackage[all,cmtip]{xy}

\usepackage{caption}

\usepackage{url}
\RequirePackage{booktabs, multirow}
\RequirePackage{pgf}
\title[A vanishing theorem for threefolds and applications]
{A vanishing theorem for threefolds in characteristic $p>5$ and applications} 
\author{Fabio Bernasconi} 
\subjclass[2020]{14E30, 14F17, 14G17, 14J17.}
\keywords{log Fano contractions, vanishing theorems, singularities, positive characteristic.}

\address{Department of Mathematics, University of Utah, Salt Lake City, UT 84112, USA} 
\email{fabio@math.utah.edu}
\newcommand{\Pic}[0]{{\operatorname{Pic}}}

\newcommand{\Diff}[0]{{\operatorname{Diff}}}

\newtheorem{thm}{Theorem}[section]

\newtheorem{cor}[thm]{Corollary}
\newtheorem{proposition}[thm]{Proposition}

\theoremstyle{definition}

\newtheorem{definition}[thm]{Definition}

\newtheorem{remark}[thm]{Remark}

\makeatletter
  
  \@addtoreset{equation}{section}
  \makeatother

\newcommand{\MO}{\mathcal{O}}

\newcommand{\Q}{\mathbb{Q}}
\newcommand{\Z}{\mathbb{Z}}

\begin{document}

\begin{abstract}
In this note we show a Kawamata-Viehweg vanishing theorem for pl-contractions on threefolds in characteristic $p>5$.
We deduce several applications for klt threefolds: the vanishing of higher direct images of structure sheaves of Mori fibre spaces, Koll\'ar's result on local vanishing and the rationality of the singularities of the base of a conic bundle.
\end{abstract}

\maketitle

%\tableofcontents

\input intro.tex
\input prelims.tex
\input vanishing.tex
\input applications.tex
\input examples.tex

\input references.tex
\end{document}

%% file: intro.tex
\section{Introduction}

The Kawamata-Viehweg vanishing theorem and its relative versions are among the most useful tools in birational geometry over fields of characteristic $0$.
One of their immediate applications is the vanishing of higher direct images for the contractions appearing in the Minimal Model Program.
However, such vanishing theorems are false in general over fields of positive characteristic as examples of log Fano contractions for which the higher direct image of the structure sheaf does not vanish have been found in \cite{Sch07, Mad16, Tan16, CT19, Ber17, Ber19}.

The aim of this note is to show that some of these pathological phenomena can happen exclusively in small characteristic for threefolds.
The main technical result is the following special case of the Kawamata-Viehweg vanishing theorem.

\begin{thm}[See Section \ref{s-kvv-pl}] \label{t-vanish-pl-contr}
	Let $k$ be a perfect field of characteristic $p>5$.
	Let $(X, \Delta)$ be a $\mathbb{Q}$-factorial dlt threefold log pair over $k$, and let $S$ be a prime Weil divisor contained in $\lfloor{\Delta \rfloor}$.
	Let $\pi \colon X \rightarrow Y$ be a projective contraction morphism between quasi-projective normal varieties over $k$ such that
	\begin{enumerate}
		\item $-(K_X+\Delta)$ is $\pi$-ample and $-S$ is $\pi$-nef;
		\item $L$ is a Weil divisor such that $L-(K_X+B)$ is $\pi$-ample for some $S \leq  B \leq \Delta$;
		\item $\dim(Y) \geq 1$.
	\end{enumerate}
	Then $R^i\pi_* \mathcal{O}_X(L) =0$ in a neighbourhood of $\pi(S)$ for $i>0$.
\end{thm}
The proof of Theorem \ref{t-vanish-pl-contr} relies on techniques developed in \cite{DH16, HW19, Ber19b} and on the Kawamata-Viehweg vanishing theorem for surfaces of del Pezzo type in characteristic $p>5$ proven in \cite{ABL20}.

Despite the various technical hypothesis, in Section \ref{s-app} we apply Theorem \ref{t-vanish-pl-contr} together with the MMP for threefolds (\cite{HX15, Bir16}) to deduce several consequences on the birational geometry of threefolds in positive characteristic.
For instance we are able to prove the vanishing of higher direct images of the structure sheaf for log Fano contractions.

\begin{thm}\label{t-main-thm}
Let $k$ be a perfect field of characteristic $p> 5$.
Let $f \colon X \rightarrow Z$ be a projective contraction morphism between quasi-projective normal varieties over $k$.
Suppose that there exists an effective $\Q$-divisor $\Delta \geq 0$ such that
\begin{enumerate}
\item $(X, \Delta)$ is a klt threefold log pair;
\item  $-(K_X+\Delta)$ is $f$-big and $f$-nef;
\item $\dim(Z) \geq 1$.
\end{enumerate} 
Then the natural map $\MO_Z \rightarrow \mathbf{R}f_*\MO_X $ is an isomorphism.
\end{thm}

We discuss two further applications of the above vanishing theorems.
In \cite{Kol11}, Koll\'ar proved a local version of the Kawamata-Viehweg vanishing theorem for klt pairs. 
We now extend his result to threefolds in characteristic $p>5$.

\begin{thm}
	Let $k$ be a perfect field of characteristic $p>5$.
	Let $(X,\Delta)$ be a threefold dlt pair over $k$.
	Let $D$ be a Weil divisor such that $D \sim_{\Q} \Delta'$, where $0 \leq \Delta' \leq \Delta$.
	Then $\mathcal{O}_X(-D)$ is Cohen-Macaulay.
\end{thm}

In Section \ref{s-examp}, we discuss singularities of the bases of Mori fibre spaces in positive characteristic.
Let us recall that over a field of characteristic $0$, the base of a klt Mori fibre space has klt singularities (see \cite[4.6]{Fu99}, \cite[Theorem 0.2]{Amb05}).
An analogous statement does not hold in positive characteristic, as examples of klt Mori fibre spaces in characteristic $2$ and $3$ such that the base has neither klt nor rational singularities have been constructed in \cite{Tan16}.
As a third application of the previous vanishing theorems, we show the following.
\begin{cor}[See Subsection \ref{ss-rat-conic}]
Let $k$ be a perfect field of characteristic $p>5$.
Let $(X, \Delta)$ be a klt threefold log pair and let $\pi \colon X \to S$ be a contraction onto a surface $S$.
If $-(K_X+\Delta)$ is $\pi$-big and $\pi$-nef, then $S$ has rational singularities.
\end{cor}

Whether the base of a threefold conic bundle has klt singularities in large characteristic still remains an open question. 
However, in Subsection \ref{ss-bad-examp} we construct examples in higher dimension of terminal Mori fibre spaces whose bases do not have log canonical singularities for every $p \geq 3$. This builds on work of Yasuda on wild quotient singularities (cf. \cite{Yas14, Yas19}). \\

\indent \textbf{Acknowledgments:} I would like to thank E. Arvidsson, C.D. Hacon, H. Tanaka, D.C. Veniani, J. Witaszek and T. Yasuda for comments and useful discussions. This work was supported by the NSF research grant n. DMS-1801851 and by a grant from the Simons Foundation; Award Number: 256202.

%% file: prelims.tex
\section{Preliminaries}

\subsection{Notation}
In this article, $k$ denotes a perfect field of characteristic $p>0$. For $m>0$, we denote by $W_n(k)$ (resp. $W(k)$) the ring of Witt vectors of length $n$ (resp. the ring of Witt vectors).

We say that $X$ is a {\em variety over} $k$ or a $k$-{\em variety} if 
$X$ is an integral scheme that is separated and of finite type over $k$. 
We denote by $D(X)$ the derived category of coherent sheaves on $X$.

We say $(X,\Delta)$ is a \emph{log pair} if $X$ is  normal variety, $\Delta$ is an effective $\Q$-divisor and $K_X+\Delta$ is $\Q$-Cartier. 
We refer to \cite{KM98, Kol13} for the basic definitions in birational geometry and of the singularities (as \emph{klt, plt}) appearing in the Minimal Model Program.
We say a morphism $\pi \colon X \to Y$ between normal varieties is a \emph{contraction} if it is proper, surjective and $\pi_*\MO_X=\MO_Y$.

\subsection{Pl-contractions}

Let us recall the definition of pl-contractions.
This type of contractions appear naturally when performing plt blow-ups (see \cite[Proposition 2.15]{GNT19}).

\begin{definition}\label{d-pl-cont}
	Let $k$ be a field.
	Let $(X, \Delta)$ be a dlt pair over $k$ and let $S$ be a prime divisor contained in $\lfloor{ \Delta \rfloor}$.
	Let $\pi \colon X \to Y$ be a projective $k$-morphism between quasi-projective normal varieties.
	We say that $\pi$ is a \emph{$(K_X+\Delta, S)$-pl-contraction} (resp. a \emph{weak $(K_X+\Delta, S)$-pl-contraction}) if
	\begin{enumerate}
		\item $-(K_X+\Delta)$ is $\pi$-ample,
		\item $-S$ is $\pi$-ample (resp. $\pi$-nef).
	\end{enumerate}
\end{definition}

In \cite[Proposition 2.8]{Ber19b}, developing ideas present in \cite{DH16}, we proved a vanishing theorem for pl-contractions whose fibers are one dimensional.
We state a slightly stronger version for our purposes.

\begin{proposition}\label{p-vanishing-small}
	Let $k$ be a perfect field of characteristic $p>5$.
	Let $(X, \Delta)$ be a $\Q$-factorial threefold dlt pair over $k$ and let $S$ be a prime divisor contained in $\lfloor{ \Delta \rfloor}$.
	Let $\pi \colon X \to Y$ be a projective contraction morphism between normal quasi-projective varieties such that
	\begin{enumerate}
		\item the maximum dimension of the fibres of $\pi$ is one;
		\item $-S$ is $\pi$-nef.
	\end{enumerate}
	Let $L$ be a Weil divisor on $X$ such that $L-(K_X+\Delta)$ is $\pi$-ample.
	Then for all $m > 0$ we have
	\[ R^1\pi_* \MO_X(L-mS)=0 \]
	in a neighbourhood of $\pi(S)$.
\end{proposition}

\begin{proof}
	The proof of \cite[Proposition 2.8]{Ber19b} applies with these hypotheses.
\end{proof}

\subsection{Rational singularities}

We recall the definition of rational singularities, following \cite{Kov17}.

\begin{definition} \label{d-rat}
	A $k$-variety $X$ has \emph{rational singularities} if
	\begin{enumerate}
		\item $X$ is normal and Cohen-Macaulay;
		\item  for every birational projective morphim $\pi \colon 
		\widetilde{X} \to X$ where $\widetilde{X}$ is Cohen-Macaulay, the natural morphism  $\MO_X \to \mathbf{R}{\pi}_* \MO_{\widetilde{X}} $ is an isomorphism.
	\end{enumerate} 
\end{definition}
In \cite[Corollary 9.11]{Kov17}, the author shows that, assuming the existence of resolution of singularities, it is equivalent to check condition (2) in Definition \ref{d-rat} in the case where $\pi$ is a resolution of singularities, thus recovering the more classical notion of rational singularities.
 
We will need the following general fact on descent of rational singularities.

\begin{proposition} \label{rationalsing}
	Let $\pi \colon Y \rightarrow X$ be a morphism of normal varieties over $k$. Assume that 
	\begin{enumerate}
		\item $Y$ has rational singularities,
		\item the natural morphism  $\mathcal{O}_X \rightarrow \mathbf{R}\pi_* \mathcal{O}_Y$ splits in $D(X)$, \label{h-2}
		\item $X$ is Cohen-Macaulay.
	\end{enumerate}
	Then $X$ has rational singularities.
\end{proposition}

\begin{proof}
	The proof closely follows \cite[Theorem 1.1]{Kov00}. 
	We prove that $X$ has pseudo-rational singularities (see \cite[Definition 1.2]{Kov17}). 
	
	Let $\varphi \colon \widetilde{X} \rightarrow X$ be projective birational morphism from a normal variety. 
	By the existence of a Macaulayfication (see \cite[Theorem 1.1]{Kaw00}) we can construct a projective birational morphism $\widetilde{Y} \to Y$ such that $\widetilde{Y}$ is Cohen-Macaulay and the following diagram commutes:
	\[
	\begin{CD}
	\widetilde{Y} @> \widetilde{\pi} >> \widetilde{X} \\
	@V \psi VV @VV \varphi V\\
	Y @> \pi >> X,
	\end{CD}
	\]
	Thus we have the following commutative diagram in $D(X)$:
		\[
	\begin{CD}
	\mathcal{O}_X @> >> \mathbf{R} \pi_* \mathcal{O}_Y \\
	@V VV @VV V\\
		\mathbf{R}\varphi_*\mathcal{O}_{\widetilde{X}} @> >> \mathbf{R} \varphi_* \mathbf{R}\widetilde{\pi}_* \mathcal{O}_{\widetilde{Y}},
	\end{CD}
	\]
	Since $Y$ has rational singularities and $\widetilde{Y}$ is Cohen-Macaulay, the composition
	\[ \mathbf{R} \pi_* \MO_Y \to \mathbf{R} \varphi_* \mathbf{R}\widetilde{\pi}_* \mathcal{O}_{\widetilde{Y}} \simeq \mathbf{R} \pi_* \mathbf{R} \psi_* \MO_{\widetilde{Y}} \simeq  \mathbf{R} \pi_* \MO_Y \]
	is a an isomorphism. 
	Therefore by assumption (\ref{h-2}) we have a splitting in $D(X)$:
	\[ \mathcal{O}_X \rightarrow \mathbf{R} \varphi_* \mathcal{O}_{\widetilde{X}} \rightarrow \mathcal{O}_X. \]
	
	Applying $\mathbf{R} \mathcal{H}om _X( - , \omega_X^\bullet)$ and Grothendieck duality to the above sequence, we have the following splitting:
	\begin{equation}\label{split-omega}
	\omega_X^\bullet \rightarrow \mathbf{R} \varphi_* \omega^{\bullet}_{\widetilde{X}} \rightarrow \omega_X^{\bullet}.
	\end{equation} 
	Since $X$ is Cohen-Macaulay, we have $ \omega_X^\bullet \simeq \omega_X[-d]. $
	Considering the $(-d)$-th cohomology group in sequence  (\ref{split-omega}),
	we deduce that the composition
	\[ \omega_X \rightarrow \varphi_* \omega_{\widetilde{X}} \rightarrow \omega_X \]
	is an isomorphism. 
	Note that $\varphi_*\omega_{\widetilde{X}}$ is a torsion-free sheaf of rank one since $\omega_{\widetilde{X}}$ is such. 
	Therefore $\varphi_* \omega_{\widetilde{X}} \rightarrow \omega_X$ is an isomorphism.
	This means that $X$ has pseudo-rational singularities. 
	By \cite[Corollary 9.11]{Kov17} we conclude that $X$ has rational singularities.
\end{proof}

\begin{remark}
	The hypothesis of $X$ being Cohen-Macaulay seems to be necessary in positive characteristic, while in characteristic $0$ it can be deduced from the other properties as an application of the Grauert-Riemenschneider vanishing theorem (see \cite[Theorem 4.3.9]{Laz04}).
\end{remark}

%% file: vanishing.tex
\section{A Kawamata-Viehweg vanishing for pl-contraction} \label{s-kvv-pl}

We now prove a special case of the Kawamata-Viehweg vanishing theorem for pl-contractions for threefolds in characteristic $p>5$. 

\begin{proof}[Proof of Theorem \ref{t-vanish-pl-contr}]
	A similar proof of \cite[Theorem 3.5]{Ber19b} applies and we show how to adapt it.
	
	Let us write $\Delta=S+\Delta'$ and $B=S+B'$, where $B' \leq \Delta'$.
	By \cite[Corollary 3.4]{HW19}, we have the following exact sequence for all $n \geq 0$:
	\[ 0 \to \mathcal{O}_X(-(n+1)S+L)\to \mathcal{O}_X(-nS+L) \to \mathcal{O}_S(G_n) \to 0, \]
	where $G_n \sim_{\Q} -nS|_S+L|_S-B_n$ for some $B_n \leq \Diff_S(B') \leq \Diff_S(\Delta')$.
	
	We now show that that $R^i \pi_* \mathcal{O}_S (G_n)=0$ for $i>0$ and for all $n>0$.
	Indeed, 
	{\small
	\[G_n \sim_{\Q} (K_{S}+ \Diff_S(B') -B_n) -(K_S+\Diff_S(B'))+ L|_S -nS|_{S} 
	\sim_{\Q}  (K_{S}+ \Gamma_n) +A_n,\] 
	\par}
	where $(S,\Gamma_n:=\Diff_S(B') -B_n)$ is dlt and $A_n:=L|_S -nS|_{S}$ is $\pi|_S$-ample.
	If $\dim(\pi(S)) \geq 1$ we conclude by \cite[Proposition 3.2]{Tan18}. If $\dim(\pi(S))=0$, then $S$ is a surface of del Pezzo type and we conclude by \cite[Theorem 1.1]{ABL20}.
	
	We now divide the proof according to the dimension of $Y$.
	
	\textbf{Case i)}. Suppose first that $\dim(Y) \geq 2$.
	To prove $R^i \pi_* \MO_X(L)=0$, it is sufficient to prove
	$R^i \pi_* \mathcal{O}_X(L-nS)=0$ for sufficiently large and divisible $n$.
	To show this, let us consider the $(-S)$-semi-ample fibration (see \cite[Lemma 2.4]{Ber19b}):
	\[\pi \colon X \xrightarrow{g} Z \xrightarrow{h} Y, \]
	and let $k>0$ be an integer such that $-kS = g^*H$, where $H$ is a $h$-ample Cartier divisor.
	By Proposition \ref{p-vanishing-small}, we have $R^ig_*\MO_X(L-nS)=0$ for $i>0$ and all $n>0$.
	By the Leray spectral sequence and the projection formula, we have:
	{\small
	\[R^i \pi_* \MO_X(L-kmS) = R^i h_* (g_* \mathcal{O}_X(L+mg^*H))= R^i h_* (g_* \mathcal{O}_X(L) \otimes \mathcal{O}_Z(mH)). \]
\par}
 Since $H$ is $h$-ample, we conclude by Serre vanishing that $R^i \pi_* \MO_X(L-kmS)=0$ for $i>0$ and for $m$ sufficiently large.
	 
	\textbf{Case ii)}. Suppose now that $\dim(Y)=1$. In this case we have that the maps $R^i \pi_* \mathcal{O}_X(L-(n+1)S) \to R^i \pi_* \mathcal{O}_X(L-nS)$ are surjective for all $n \geq 0$.
	Let $z=\pi(S)$ and let $m>0$ such that $\pi^*z=mS$.
	In particular the following map
	\[R^i \pi_* \MO_{X}(L-mS) \simeq R^i \pi_* \MO_X(L) \otimes \MO_X(-z) \to R^i \pi_* \MO_X(L) \]
	is surjective.
	We conclude that $R^i \pi_* \MO_X(L)$ is zero in a neighbourhood of $z$ by Nakayama's lemma.
\end{proof}

%% file: applications.tex
\section{Applications}\label{s-app}

In this section we describe some applications of Theorem \ref{t-vanish-pl-contr} to the birational geometry of threefolds in positive characteristic.

\subsection{Vanishing for log Fano contractions}

We prove the vanishing of higher direct images of the structure sheaf for log Fano contractions in characteristic $p>5$.

\begin{proof}[Proof of Theorem 1.2]
	It is sufficient to prove $R^if_*\MO_X=0$ for $i>0$.
	By \cite[Theorem 3.3]{Tan18} and \cite[Corollary 1.8]{BT19}, there exists an open subset $U \subset Z$ such that $(R^i f_* \MO_X)|_{U}=0$ and $Z \setminus U$ is a finite set of points.
	Let $z$ be a closed point in $Z \setminus U$.
	By \cite[Proposition 2.15]{GNT19},
	there exists a birational morphism
	$\pi \colon Y \to Z$
	and an effective $\Q$-divisor $\Delta_Y$ on $Y$ such that 
	\begin{enumerate}
		\item[(i)] $(Y, \Delta_Y)$ is a $\Q$-factorial plt pair (in particular, $Y$ is klt),
		\item[(ii)] $S:=(\pi^{-1}(z))_{\text{red}}$ is an irreducible component of $\lfloor{ \Delta_Y \rfloor}$,
		\item[(iii)] $\pi$ is a weak $(K_Y+\Delta_Y,S)$-pl-contraction.
	\end{enumerate}
	Let us consider the following diagram
	\[
	\begin{CD}
	W @> \varphi >> Y \\
	@V \psi VV @VV \pi V\\
	X @> f >> Z,
	\end{CD}
	\]
	where $\varphi$ and $\psi$ are log resolutions.
	Since $X$ and $Y$  have rational singularities by \cite[Corollary 1.3]{ABL20}, to prove $R^i f_* \MO_X=0$ in a neighbourhood of $z$ it is sufficient to prove $R^i \pi_* \MO_Y=0$ for $i>0$. 
	This follows from Theorem \ref{t-vanish-pl-contr}.
\end{proof}

As an application we discuss lifting of log Fano contractions in mixed characteristic.
Let us start by recalling the following lifting result:

\begin{thm}\label{lift-vanishing}
Let $k$ be a pefect field and
let $f \colon X \to Y$ be a morphism of $k$-schemes such that $\mathbf{R}f_* \MO_X \simeq \MO_Y$.
If $X$ lifts to $W_m(k)$ (resp. formally lifts to $W(k)$), then the morphism $f \colon X \to Y$ lifts to $W_m(k)$ (resp. formally lifts to $W(k)$).
\end{thm}
\begin{proof}
See \cite[Theorem 3.1]{CvS09}.
\end{proof}

This result has been used by Hacon and Witaszek to show the lifting of birational contractions between klt pairs (see \cite[Corollary 5.1]{HW19}).
Combining Theorem \ref{lift-vanishing} with Theorem \ref{t-main-thm}, we extend their result to show the lifting of log Fano contractions.

\begin{cor}\label{lifting}
Let $k$ be a perfect field of characteristic $p>5$.
Let $(X, \Delta)$ be a klt threefold log pair and let $f \colon X \to Z$ be a projective contraction morphism over $k$ between quasi-projective varieties.
Assume that
\begin{enumerate}
\item $X$ lifts to $W_m(k)$ (resp. formally lifts to $W(k)$),
\item $-(K_X+\Delta)$ is $f$-big and $f$-nef,
\end{enumerate} 
Then the morphism $f \colon X \to Z$ lifts to $W_m(k)$ (resp. formally to $W(k)$).
\end{cor}

\subsection{Depth of sheaves on klt pairs}

In \cite[Theorem 2]{Kol11}, Koll\'ar showed that certain divisorial sheaves on klt pairs over a field of characteristic $0$ are Cohen-Macaulay, and he interprets this result as a local Kawamata-Viehweg vanishing theorem.
In positive characteristic, a version for strongly $F$-regular pairs has been proved by Patakfalvi and Schwede (see \cite[Theorem 3.1]{PS14}).

As an application of Theorem \ref{t-vanish-pl-contr}, we extend Koll\'ar's result to the case of klt threefolds in characteristic $p>5$.

\begin{thm}\label{t-local-kvv}
	Let $k$ be a perfect field of characteristic $p>5$.
	Let $(X,\Delta)$ be a threefold dlt pair over $k$.
	Let $D$ be a $\Z$-divisor such that $D \sim_{\Q} \Delta'$, where $0\leq \Delta' \leq \Delta$.
	Then $\mathcal{O}_X(-D)$ is Cohen-Macaulay.
\end{thm}

\begin{proof}
	If $X$ is $\Q$-factorial, then we conclude by \cite[Corollary 1.3]{ABL20}.
	We can thus work on an affine neighbourhood of a non $\Q$-factorial closed point $x$.
	Using \cite[Proposition 2.43]{KM98} and resolution of singularities for threefolds in positive characteristic, we may assume $(X, \Delta)$ is klt. Let us write $\Delta=\Delta'+\Delta''$.
	Let $\pi \colon Y \to X$ be a plt blow-up at $x$ (see \cite[Proposition 2.15]{GNT19}), i.e.:
	\begin{enumerate}
		\item[(i)] $S=(\pi^{-1}(x))_{\text{red}}$ is the unique $\pi$-exceptional divisor,
		\item[(ii)] $(Y, S+\pi_*^{-1}\Delta)$ is a $\Q$-factorial plt pair,
		\item[(iii)] $\pi$ is a weak $(K_Y+ S+\pi_*^{-1}\Delta, S)$ pl-contraction.
	\end{enumerate}
	Consider $\pi^*(D-\Delta')=\pi_*^{-1}D - \pi_*^{-1}\Delta' +cS$ for some $c \in \Q$.
	
	Define $D_Y:=\pi_*^{-1}D - \lfloor cS \rfloor \sim_{\Q} \pi_*^{-1}\Delta' + \{cS \} + \pi^*(D-\Delta')$.
	Note that $K_Y+D_Y-(K_Y+S+\pi_*^{-1}\Delta') \sim_{\Q} \{ cS \}-S+\pi_*^{-1}(D-\Delta')$ is nef and big over $X$, so $R^i \pi_*\mathcal{O}_Y(K_Y+D_Y)=0$ for $i>0$ by Theorem \ref{t-vanish-pl-contr}.
	
	We have
	\[K_Y+\pi_*^{-1}\Delta=\pi^*(K_X+\Delta) +aS,\]
	for some $a>-1$.
	We write $aS+\left\{ cS \right\}=\lceil{aS+\left\{ cS \right\} \rceil}-A_Y$.
	We have the following equalities:
	{\small
		\begin{align*}
			\lceil{aS+\left\{ cS \right\} \rceil}-D_Y & \sim_{\Q} aS+ \left\{ cS \right\} +A_Y - (\pi_*^{-1}\Delta' + \{cS \} + \pi^*(D-\Delta')) \\
			& \sim_{\Q} K_Y +\pi_*^{-1}\Delta +A_Y - \pi_*^{-1}(\Delta')-\pi^*(K_X+\Delta+D-\Delta')\\
			& \sim_{\Q} K_Y +\pi_*^{-1}\Delta^{''} +S + (A_Y-S) -\pi^*(K_X+\Delta''+D).
		\end{align*}
		\par}
	Thus we conclude again by Theorem \ref{t-vanish-pl-contr} that $R^i \pi_* \mathcal{O}_Y(	\lceil{aS+\left\{ cS \right\} \rceil}-D_Y )=0$.
	Moreover, by Fujita's lemma (see \cite[Lemma 7.30]{Kol13}), we have $\mathcal{O}_X(-D)=\pi_* \mathcal{O}_Y(-D_Y)=\pi_* \mathcal{O}_Y(\lceil{aS+\left\{ cS \right\} \rceil}-D_Y)$.
	
	Let us sum up the vanishing results just proven:
	\begin{enumerate}
		\item The support of $\mathcal{O}_Y(-D_Y)$ is $Y$ and $\mathcal{O}_Y(-D_Y)$ is Cohen-Macaulay by \cite[Corollary 1.2]{ABL20} since $Y$ is klt and $\Q$-factorial;
		\item we have $\mathcal{H}om(\mathcal{O}_Y(-D_Y), \omega_Y)\simeq \mathcal{O}_Y(K_Y+D_Y)$ since both are reflexive sheaves on a normal variety which coincide on a big open set. Thus $R^i\pi_* \mathcal{H}om(\mathcal{O}_Y(-D_Y), \omega_Y)=0$ for $i>0$;
		\item the composition in the derived category
		\[\qquad \pi_* \mathcal{O}_Y(-D_Y) \to \mathbf{R}\pi_*\mathcal{O}_Y(-D_Y) \to \mathbf{R}\pi_*\mathcal{O}_Y(\lceil{aS+\left\{ cS \right\} \rceil}-D_Y)\] 
		is an isomorphism.
	\end{enumerate}
	Thus by \cite[Theorem 2.74]{Kol13}, we conclude that $\pi_* \mathcal{O}_Y(-D_Y)=\mathcal{O}_X(-D)$ is Cohen-Macaulay.
\end{proof}

\begin{remark}
	The assumption on the characteristic in Theorem \ref{t-local-kvv} is optimal as for $p \leq 5$ there are examples of klt threefold singularities which are not Cohen-Macaulay (see \cite{CT19,Ber17,ABL20}).
\end{remark}

%% file: examples.tex
\section{Singularities of the bases of Mori fibre spaces}\label{s-examp}

In this section we discuss singularities of the base of Mori fibre spaces in positive characteristic. 

\subsection{Threefold conic bundles} \label{ss-rat-conic}

In \cite[Theorem 3.8]{NT}, the authors prove that the base scheme of a threefold conic bundle has $W\MO$-rational singularities over perfect fields of characteristic $p > 5$.
Using the recent result of Hacon and Witaszek on the validity of the MMP in characteristic $5$ (\cite{HW19}), this can extended to characteristic five.

We present a stronger result: the base of a threefold conic bundle has rational singularities for $p>5$.

\begin{thm}
	Let $k$ be a perfect field of characteristic $p>5$.
	Let $(X, \Delta)$ be a klt threefold log pair and let $\pi \colon X \to S$ be a projective contraction onto a surface $S$.
	If $-(K_X+\Delta)$ is $\pi$-big and $\pi$-nef, then $S$ has rational singularities. 
\end{thm}

\begin{proof}
	By \cite[Corollary 1.3]{ABL20}, $X$ has rational singularities and $\MO_S \to \textbf{R}\pi_*\MO_X$ is a quasi-isomorphism by Theorem \ref{t-main-thm}.
	Since $S$ is a normal surface, it is in Cohen-Macaulay and thus we conclude by Proposition \ref{rationalsing}.
\end{proof}

\begin{remark}
	It would be natural to expect that the base of a Mori fibre space has klt singularities, at least for sufficiently large $p>0$.
	However, even in the case of threefold conic bundles, we have  little evidence.
	Let us note that Koll\'ar proved that the base is smooth if the total space is smooth (see \cite[Complement 4.11.2]{Kol91}).
\end{remark}

\subsection{Pathological examples in higher dimension} \label{ss-bad-examp}

The aim of this section is to show new examples of Mori fibre spaces whose bases have bad singularities in characteristic $p>0$.
The first examples were constructed for $p=2,3$  in \cite[Theorem 1.1]{Tan16}.

Our construction is based on the work of Yasuda on wild quotient singularities (see \cite{Yas14, Yas19}).
We fix an algebraically closed field $k$ of characteristic $p>0$ and the cyclic group $G:= \mathbb{Z}/p\mathbb{Z}$. Let us recall that for every integer $1 \leq i \leq p$ we have a unique indecomposable representation of $G$ over $k$ on a $k$-vector space of dimension $i$ denoted by $V_i$. Let $V$ be a $G$-representation and consider the quotient $V \rightarrow X:=V/G$. The representation $V$ decomposes into sum of indecomposable ones
\[V= \bigoplus_{\lambda=1}^l V_{d_\lambda}\] with $1 \leq d_\lambda \leq p$ and we introduce the following invariant for $X$:
\[D_V = \sum_{i=1}^l \frac{d_\lambda(d_{\lambda}-1)}{2}. \]
\begin{thm}[{\cite{Yas19}}]\label{yasuda}
Suppose $D_V \geq 2$. Then the quotient variety $X$ is terminal (resp. canonical, log canonical) if and only if $D_V >p$ (resp. $D_V \geq p$, $D_V \geq p-1$). 
\end{thm} 
We now construct examples of Mori fibre spaces whose bases do not even have log canonical singularities.
\begin{thm}\label{baseNotLC}
Let $k$ be an algebraically closed field of characteristic $p \geq 5$. Then there exists a projective contraction $f \colon X \rightarrow Y$ of normal $k$-varieties such that
\begin{enumerate}
\item $X$ is a $\Q$-factorial terminal quasi-projective variety of dimension $p+3$; 
\item $Y$ is a $\Q$-factorial affine variety of dimension three which is not log canonical;
\item $\rho(X/Y)=1$ and $-K_X$ is $f$-ample, equivalently $f$ is a Mori fibre space.
\end{enumerate}
\end{thm}

\begin{proof}
Consider the indecomposable representation of $G$ on the three-dimensional space $V_3$. Since $D_{V_3} = 3 < p-1$ we have by Theorem \ref{yasuda} that the quotient $Y:= V_3/ G$ is not log canonical. \\
\indent Now we consider the $(\mathbb{P}^1_k)^{p}$ with the following $G$-action:
{\small
\[T( ([x_1: y_1], [x_2: y_2], \dots , [x_p:y_p])) = ([x_p:y_p], [x_{1}: y_{1}], \dots, [x_{p-1}: y_{p-1}]). \]
\par}
Let $G$ act diagonally on $(\mathbb{P}^1_k)^{p} \times V_3$. 
In this case the quotient $X:=((\mathbb{P}^1)^{p} \times V_3)/G$ has terminal singularities. 
Indeed, on a local chart $\mathbb{A}^p \subset (\mathbb{P}^1)^p$ the action is the sum of the irreducible representations $V_p \oplus V_3$ and since $D_{V_{p} \oplus V_3} = \frac{p(p-1)}{2} + 3 \geq p+1$ we conclude $X$ is terminal by Theorem \ref{yasuda}. 
Consider now the following diagram
\begin{equation*}
\xymatrix{
(\mathbb{P}^1)^{p} \times V_3 \ar[r]^{ \qquad \pi} \ar[d] & X \ar[d]^{f} \\
V_3 \ar[r] & Y
}
\end{equation*}
Since $\pi$ is \'etale in codimension one, we deduce $-K_X$ is $f$-ample and we are only left to prove that the relative Picard number $\rho(X/Y)$ is one.
This is immediate since $\Pic(X) \hookrightarrow (\Pic ((\mathbb{P}^1)^{p} \times V_3))^{G} \simeq \mathbb{Z}.$
\end{proof}

\begin{remark}
The totale space appearing in the previous construction do not have Cohen-Macaulay singularities by \cite{ES80}.
Let us note moreover that in Theorem \ref{baseNotLC}, the relative dimension of the fibration increases with the characteristic $p$.
\end{remark}